\documentclass[11pt]{article}
\makeindex
\usepackage{graphicx}
\usepackage[latin1]{inputenc} \usepackage[T1]{fontenc} \usepackage{lmodern}
\usepackage{pgf}
\usepackage{epsfig,graphicx,color,amsmath,a4wide,amssymb,amsfonts,amsbsy,xy,latexsym,epsf,pstricks,verbatim,times}

\usepackage{pgf,color}
\usepackage{pgfarrows}

\usepackage{changebar}

\definecolor{LightGrey}{rgb}{.85,.85,.85}
\definecolor{DarkGrey}{rgb}{.5,.5,.5}

\definecolor{Blue}{rgb}{.0,.0,0.9}
\definecolor{LightBlue1}{rgb}{.2,.4,0.9}
\definecolor{LightBlue2}{rgb}{.3,.5,0.9}
\definecolor{LightBlue3}{rgb}{.4,.6,0.9}
\definecolor{LightBlue4}{rgb}{.5,.7,.9}
\definecolor{LightBlue5}{rgb}{.6,.8,.9}
\definecolor{LightBlue6}{rgb}{.7,.9,.9}

\definecolor{Red}{rgb}{.9,.0,.0}
\definecolor{LightRed1}{rgb}{0.9,.2,.4}
\definecolor{LightRed2}{rgb}{0.9,.3,.5}
\definecolor{LightRed3}{rgb}{0.9,.4,.6}
\definecolor{LightRed4}{rgb}{.9,.5,.7}
\definecolor{LightRed5}{rgb}{.9,.6,.8}
\definecolor{LightRed6}{rgb}{.9,.7,.9}

\def\NS{\mathop{\rm NS }}

\definecolor{Grey}{rgb}{.5,.5,.5}

\definecolor{Blue}{rgb}{.0,.0,0.9}
\definecolor{LightBlue1}{rgb}{.2,.4,0.9}
\definecolor{LightBlue2}{rgb}{.3,.5,0.9}
\definecolor{LightBlue3}{rgb}{.4,.6,0.9}
\definecolor{LightBlue4}{rgb}{.5,.7,.9}
\definecolor{LightBlue5}{rgb}{.6,.8,.9}
\definecolor{LightBlue6}{rgb}{.7,.9,.9}

\definecolor{Red}{rgb}{.9,.0,.0}
\definecolor{LightRed1}{rgb}{0.9,.2,.4}
\definecolor{LightRed2}{rgb}{0.9,.3,.5}
\definecolor{LightRed3}{rgb}{0.9,.4,.6}
\definecolor{LightRed4}{rgb}{.9,.5,.7}
\definecolor{LightRed5}{rgb}{.9,.6,.8}
\definecolor{LightRed6}{rgb}{.9,.7,.9}

\xyoption{all}
\usepackage[english]{babel}

\newcounter{noalgo}[section]
\newdimen\indentalgo
\newdimen\indentalgodec\indentalgo=0.0mm\indentalgodec=10mm

\newcommand{\If}{\advance\indentalgo by \indentalgodec {\bf if }}

\newcommand{\For}{\global\advance\indentalgo by \indentalgodec {\bf for }}

\newcommand{\Endindent}{\global\advance\indentalgo by -\indentalgodec}

\newdimen\decalage \decalage=0.5cm
\newcounter{algo} 

\newcommand{\PP}{\mathbb P}
\newcommand{\hPP}{\hat {\mathbb P}}

\def\<<{\leavevmode
  \raise0.28ex\hbox{$\scriptscriptstyle\langle\!\langle$}\nobreak
  \hskip -.6pt plus.3pt minus.2pt\,}
\def\>>{\,\nobreak\hskip -.6pt plus.3pt minus.2pt
  \raise0.28ex\hbox{$\scriptscriptstyle\rangle\!\rangle$}}

\def\Proj{\mathop{\rm{Proj}}\nolimits }

\def\cS{{\cal S}}

\newtheorem{lemma}{Lemma}

\providecommand{\myproofname}{Proof}

\usepackage{enumitem}

\begin{document}

%\author{Jean-Marc Couveignes
%D{\'e}partement de Math{\'e}matiques et Informatique,
%Universit{\'e} Toulouse 2, 5 all{\'e}es Antonio Machado, 31058
%Toulouse     c{\'e}dex 9.} 
\author{Jean-Marc Couveignes\thanks{Corresponding author, Universit\'e de     Toulouse le
    Mirail,
5 allées Antonio Machado, 31058 Toulouse Cédex 9, France.  email :
couveig@univ-tlse2.fr, telephone 33561504131, fax 33561504173.}
\thanks{INRIA Bordeaux Sud-Ouest, Institut de Mathématiques de Bordeaux, CNRS.}
%D{\'e}partement de Math{\'e}matiques et Informatique,
%Universit{\'e} Toulouse 2, 5 all{\'e}es Antonio Machado, 31058
%Toulouse     c{\'e}dex 9.} 
and Jean-Gabriel Kammerer\thanks{\textsc{DGA}/\textsc{MI}, BP
    7 35998 RENNES ARMEES.}~\thanks{IRMAR, Universit\'e de Rennes 1, Campus de
    Beaulieu, F-35042 Rennes.}}

\title{The Geometry of Flex
Tangents to a  Cubic Curve and 
 its Parameterizations\thanks{Research supported by 
the ``Agence Nationale de la Recherche''  through project  ALGOL (ANR-07-BLAN-0248) and by
the French ``D{\'e}l{\'e}gation
    G{\'e}n{\'e}rale pour l'Armement''.}}

\maketitle

\bibliographystyle{plain}

\begin{abstract}
We show how the study of the geometry of the nine
flex tangents to  a cubic produces
pseudo-parameterizations, including the ones given
by Icart, Kammerer, Lercier, Renault and Farashahi, and infinitely many
new ones.
\end{abstract}
\hfill {\it To Jean-Jacques Quisquater, on the occasion of his {\'e}m{\'e}ritat}
\section{Introduction}\label{sec:intro}
Much attention has been focused recently
on the problem of computing
points on a given elliptic curve over a finite
field  in deterministic polynomial time. 
This problem arises in a very natural manner in many cryptographic protocols
when one wants to encode messages into the group of points of an elliptic
curve. A good example of the algorithmic and cryptologic motivations in finding
these parameterizations can be found in the identity-based encryption from
\cite{IDEnc}. The difficulty is to deterministically find a field element
$x$ such that some polynomial in $x$ is a square, see \cite{Koblitz}, Section~6.1.8. For example, when the curve is given by a reduced Weierstrass equation
$y^2 = x^3 + ax +b$, we deterministically search $x$ such that $x^3+ax+b$ is a square in the
field. 

In 2006, Shallue and Woestjine \cite{SW} proposed 
a first practical deterministic algorithm.
In 2009, Icart \cite{icart} proposed  another
deterministic encoding for elliptic curves over
a field $k$ with $q$ elements,  when $q$ is congruent to $2$ modulo $3$.
Icart's algorithm has quasi-quadratic   complexity 
in $\log q$.
Kammerer, Lercier and Renault \cite{KLR} proposed a different 
encoding 
under the additional condition that the elliptic curve has a rational point of
order $3$, and even for a special class of hyperelliptic curves. Farashahi
\cite{fara} found yet another parameterization for such elliptic curves too. 
A crucial  point in \cite{icart,KLR,fara} is that the map 
$x\mapsto x^3$ is bijective for
a finite field $k$ having  cardinality congruent to $2$ modulo $3$.
Its inverse map is $x\mapsto x^{e}$ where $e\bmod q-1$ is the 
inverse of $3\bmod q-1$  and $0\le e<q-1$.
Exponentiation by $e$ can be computed in deterministic
time $(\log q)^{2+o(1)}$ using the 
fast exponentiation algorithm.   So in order to deterministically 
compute points on an elliptic curve $C$ over such a finite field, one can afford the usual
field operations together with   cubic roots.
In other words, one looks for a
parameterization of the elliptic curve
by cubic radicals. Such a parameterization will be called
a {\it pseudo-parameterization} in this article.
Finding such a pseudo-parameterization is a special case of the
problem of finding parameterizations of curves by radicals \cite{SWP}.

We show how  such pseudo-parameterizations can be obtained from the
study of the dual curve of the elliptic curve $C$. 
In a nutshell, we produce points on $C$ as intersection points
between $C$ and well chosen lines. If $D$ is a line in the 
projective plane, then the intersection $D.C$ consists of three
points, counting multiplicities.
These three points can be computed  by solving a cubic equation.
We recall in Section~\ref{sec:cubiceq}  how to 
derive the  Tartaglia-Cardan formulae for this purpose.
Recall  these formulae run in two steps.
One first has to compute a square root of the discriminant. The three
solutions 
are then calculated using the field operations and cubic roots.
Since  cubic roots are not a problem in our context, the only
remaining difficulty  is computing the square root of  the discriminant. So we  choose the 
line $D$ in such a way that the discriminant of the intersection
$D.C$ is a square, and we assume that we have an algebraic  formula for its
square root. More precisely, we consider a line $D_t$ depending
on a rational formal parameter $t$. This means that
the  coefficients in the
projective equation of $D_t$ are polynomials in the indeterminate
$t$. The discriminant $\Delta(t)$ of the intersection $L_t.C$ is then a rational
fraction in $t$. We ask that this discriminant be a square in $k(t)$.
We  compute once for all a  formal square root $\delta(t)$
of $\Delta(t)$. For every value of $t$ we can then produce a point
on $C$ using only the field operations and cubic roots.

We recall in Section~\ref{section:dual} that the  projective lines in $\PP$ are parametrized 
by the dual plane
$\hPP$. The line in $\PP$ with projective equation $UX+VY+WZ=0$ is represented
by the point $[U:V:W]\in \hPP$.
A rational family of lines $t\mapsto D_t$ thus gives  rise
to a rational curve $L$ inside $\hPP$. Indeed, if the projective
equation  of $D_t$ is $U(t)X+V(t)Y+W(t)Z=0$ then the map
$t\mapsto [U(t):V(t):W(t)]$ parametrizes a rational curve
inside $\hPP$.
The discriminant $\Delta(t)$ vanishes whenever $D_t.C$ has a multiple
root. This happens if and only if $D_t$ is {\it tangent} to $C$. Not
every projective line is tangent to $C$. The subset of $\hPP$ corresponding
to lines that are tangent to $C$ is a curve denoted $\hat C$ and called
the {\it dual curve of $C$}. So  $\Delta(t)$ describes the intersection between
the rational curve $L$ and the dual curve $\hat C$.
And $\Delta(t)$ is a square if and only if every
 point in the intersection between $L$ and $\hat C$ has even
multipicity. So we will be interested in rational curves 
$L$ in $\hPP$  that have even intersection with the dual
curve to the cubic curve $C$.  The connection between such 
curves and pseudo-parameterizations is detailed in Section~\ref{sec:parame}.

Because the dual curve $\hat C$  plays such an important role we will
study it in Section~\ref{section:dual}. This curve has genus $1$ and
$9$ singularities, all cusps. Indeed the nine cusps of $\hat C$ correspond
to the nine flex tangents to $C$, while the smooth points on $\hat C$
parametrize  the tangent lines to $C$ that are not flexes. These nine 
points in the dual plane form an interresting configuration that we study
in Section~\ref{sec:flexes}. We are particularly interested  in 
rational curves $L$   passing through several
among these nine points. We will find that  many
such  curves $L$ have even intersection with $\hat C$.
We  will show in Section~\ref{sec:inter} that these curves
give rise to all the known pseudo-parameterizations  of $C$ 
found by Icart, Farashahi,  Kammerer, Lercier,
Renault,
and to 
several new ones.  It is then  natural
to ask how many rational curves on $\hPP$ have even intersection
with $\hat C$. 
We shall see in Section~\ref{sec:K3} that there
are infinitely many such rational curves, giving rise
to infinitely many inequivalent pseudo-parameterizations.
These curves 
lift  to rational curves on the degree
two covering $\Sigma$ of the dual plane ramified along  $\hat C$.
This will lead us to the classical and beautiful topic of 
rational curves on $K3$ surfaces.

Throughout the paper, we denote by $k$ a field with characteristic different
from $2$ and $3$, by $\bar k \supset  k$ an algebraic closure of $k$, and by
$\zeta_3\in \bar k$ a primitive third root of unity.
We set $\sqrt{-3}=2\zeta_3+1$.

The Maple \cite{Wa} code for the calculations in this article can be found on the
authors'  web pages.

\section{Solving cubic equations}\label{sec:cubiceq}

In this section we recall the Tartaglia-Cardan formulae for solving
cubic
equations by radicals. A modern treatment can be found in \cite{DF}.
We believe it is worth stating these equations
in an unambiguous form, that is well adapted to our context, and does
not make excessive use of radicals and roots of unity. In other words
we need regular and generic formulae.
Let $h(x)=x^3-s_1x^2+s_2x-s_3$ be a degree $3$ separable polynomial
in $k[x]$. Call $r_0$, $r_1$ and $r_2$ the three roots of $h(x)$ in 
$\bar k$. Set \[\delta=\sqrt{-3}(r_1-r_0)(r_2-r_1)(r_0-r_2)\] and 
$\Delta =\delta^2$. Note that $\Delta$ is 
the usual discriminant multiplied by $-3$. We call it
the {\it twisted discriminant}.
Since it is a symmetric function
of the roots, it
 can be expressed as a polynomial in
$s_1$, $s_2$ and $s_3$. Indeed 
\[\Delta = 81s_3^2-54s_3s_1s_2-3s_1^2s_2^2+12s_1^3s_3+12s_2^3.\] In particular $\Delta$
lies in $k$. 
Let  $l=k(\zeta_3,\delta)\subset \bar k$ be the field obtained by
adjoining $\delta$ and a primitive third root of unity to $k$.
We set $m=l(r_1,r_2,r_0)$.

If the extension $l\subset m$ is non-trivial then it is
a cyclic
cubic extension. 
Since $l$ contains a primitive third root of unity,
this cubic extension is a Kummer extension: it is generated by the cubic
root of some element in $l$. 
Let $\sigma$
be the generator of the Galois group 
that sends $r_i$ to $r_{i+1}$ for
$i\in \{0,1,2\}$, with the convention that indices make sense modulo $3$.
We set \[\rho = r_0+\zeta_3^{-1}r_1+\zeta_3^{-2}r_2\]
and we check  that $\sigma(\rho)=\zeta_3\rho$.
We set $R=\rho^3$ and we check that $R$ is invariant by $\sigma$.
So $R$ is an invariant for the alternate group 
acting on $\{r_1,r_2,r_3\}$ and
it can be expressed as a polynomial in $s_1$, $s_2$, $s_3$ and
$\delta$. 
Indeed we find
\[R =\rho^3= s_1^3+\frac{27}{2}s_3-\frac{9}{2}s_1s_2-\frac{3}{2}\delta.\]
Similarly we set \[\rho' = r_0+\zeta_3r_1+\zeta_3^{2}r_2\]
and we check that 
\[R' = \rho'^3  = s_1^3+\frac{27}{2}s_3-\frac{9}{2}s_1s_2+\frac{3}{2}\delta.\]

We note  that
$\rho\rho'=r_0^2+r_1^2+r_2^2-r_0r_1-r_1r_2-r_2r_0$ is invariant
by the full symmetric group and is indeed equal to $s_1^2-3s_2$.
So both $\rho$ and $\rho'$ are computed by extracting a single 
cubic root.

Finally, the three roots $r_0$, $r_1$, $r_2$ can be expressed
in terms of $\rho$ by solving the linear system:
\[\left\{
\begin{array}{ccc}
r_0+r_1+r_2 &=&s_1\\
r_0+\zeta_3^{-1}r_1+\zeta_3r_2 &=&\rho\\
r_0+\zeta_3r_1+\zeta_3^{-1}r_2 &=& \rho'
\end{array}
\right.\]
In particular the  formula for the root 
\begin{equation}\label{eq:r0}
r_0=\frac{s_1+\rho+\rho'}{3}
\end{equation}
 does
not involve $\zeta_3$.

\section{The dual curve of a cubic}\label{section:dual}

In this section we review the properties of the dual of a 
cubic curve. A thorough treatment of the duality for plane
curves can be found in \cite{Ge-Ka-Ze} , \cite{HKT} and  \cite{H}.
Let $E=k^3$ and let
$\hat E$ be the dual of $E$.
Let $U=(1,0,0)$, $V=(0,1,0)$ and $W=(0,0,1)$. So
$(U,V,W)$ is the  
canonical basis of $E$.
Let $(X,Y,Z)$ be the 
dual basis of $(U,V,W)$.
Let $\PP=\Proj(E)=\Proj k[X,Y,Z]$ be the projective plane
over $k$. Let $\hat \PP=\Proj(\hat E)= \Proj k[U,V,W]$ be the dual
projective plane. 
The main idea of projective dualy is that
points in $\hPP$ parametrize lines in $\PP$, and conversely.
The point $[U:V:W]$ in $\hat \PP$ corresponds to the
line with equation $UX+VY+WZ=0$ in $\PP$. And the point
$[X:Y:Z]$ in $\PP$ parametrizes the line $XU+YV+ZW=0$ in $\hPP$.

Now let $C\subset \PP$ be an absolutly integral  curve with equation
$F(X,Y,Z)=0$. Let
$F_X=\frac{\partial F}{\partial X}$,
$F_Y=\frac{\partial F}{\partial Y}$,
$F_Z=\frac{\partial F}{\partial Z}$ be the three partial derivatives
of $F$.
The tangent to $C$ at a smooth  point $P=[X_P:Y_P,Z_P]$ has equation
\[F_X(X_P,Y_P,Z_P)U+F_Y(X_P,Y_P,Z_P)V+F_Z(X_P,Y_P,Z_P)W=0.\]
The corresponding point in $\hPP$ is 
$[F_X(X_P,Y_P,Z_P):F_Y(X_P,Y_P,Z_P):F_Z(X_P,Y_P,Z_P)]$.
The Zariski closure of the set of all such points is
the {\it dual}  $\hat C$ of $C$. So 
 $\hat C$ is the closure of the image of the so called Gauss morphism
\[
\xymatrix{
\omega_C : & C^{smo} \ar@{->}[r]  & \hat \PP\\
    & [X:Y:Z]\ar@{|->}[r]&[F_X(X,Y,Z),F_Y(X,Y,Z),F_Z(X,Y,Z)],
}
\]
where $C^{smo}$ is the locus of smooth points on $C$.

We assume that the characteristic of  $k$ is odd, and that not
every point on the curve $C$ is a flex or a singular point
(in particular $C$
is not a line). Then $\hat C$ is an absolutely
integral curve.
And the dual of $\hat C$ is $C$. This is the biduality  theorem \cite[Theorem 5.91]{HKT}. Duality is very useful because it translates properties
of $C$ into properties of $\hat C$ and conversely. In particular
the Gauss map $\omega_C$ is a birational map from  $C$
to  $\hat C$. It maps the  flexes of $C$ onto 
the cusps of $\hat C$.

The first non-trivial example of duality concerns conics
(smooth 
plane projective curves of degree $2$). The
dual of  conic is a conic.

We now assume that $C$ is a smooth cubic.
Then  $\hat C$  has degree $6$ and to 
 each of the nine flexes of $C$ there corresponds
an ordinary cusp on $\hat C$. Since $\hat C$ has geometric genus
$1$ and arithmetic genus $10=(6-1)(6-2)/2$ we deduce that there is no other
singularity on it than these nine cusps.
For example, if $C$  has  equation $F(X,Y,Z)=0$ where 
\begin{equation}\label{eq:Hessian}
F(X,Y,Z)=X^3+Y^3+Z^3-3aXYZ,
\end{equation}
then the dual curve has equation $G(U,V,W)=0$ where 
\begin{equation}\label{eq:G}
\begin{split}
G(U,V,W)=U^6+V^6+W^6
-6a^2(U^4VW+UV^4W+UVW^4)\\+(4a^3-2)(U^3V^3+U^3W^3+V^3W^3)+(12a-3a^4)U^2V^2W^2.
\end{split}
\end{equation}

\begin{figure}[h]
\centering
{\includegraphics[width=4cm,height=4cm]{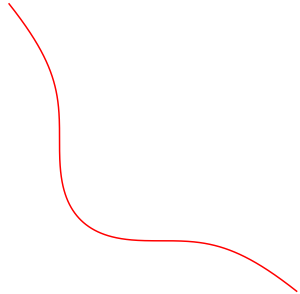}}
\caption{The cubic with equation $X^3+Y^3+Z^3=0$}
\label{fig:cub}
\end{figure}

\begin{figure}[h]
\centering
{\includegraphics[width=4cm,height=4cm]{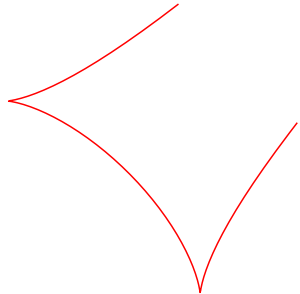}}
\caption{The dual curve with equation $ U^6+V^6+W^6-2U^3V^3-2V^3W^3-2U^3W^3 =0$}
\label{fig:dual}
\end{figure}

The equation of the dual is found by eliminating $X$, $Y$, and $Z$
in the system
\[\left\{
\begin{array}{ccc}
U &=&F_X(X,Y,Z)\\
V &=&F_Y(X,Y,Z)\\
W &=&F_Z(X,Y,Z)
\end{array}
\right.\]
%To this end, we may use  the Maple code in Figure~\ref{fig:MapleI}.
The real loci of the two
curves $C$  and $\hat C$ are represented
in Figure~\ref{fig:cub} and Figure~\ref{fig:dual}  respectively in
the case $a=0$.

The equation of the dual curve  arises naturally when one studies
the intersection of the cubic $C$ with a projective line $D$.
Indeed such a line $D\subset \PP$ meets $C$ in exactly three points unless
it is a tangent line to $C$ (in which case we have one  simple point
and one double point) or even a flex (in which case we have one 
triple point).
Assume that  $D$ is the line with equation
\begin{equation}\label{eq:line}
UX+VY+WZ=0.
\end{equation}
The intersection   
$D.C$ is described by the homogeneous system consisting of
Equation~(\ref{eq:line}) and the equation of the cubic $C$.
We can use Equation~(\ref{eq:line}) to eliminate one of the three
variables $X$, $Y$, $Z$ in  the equation
of  $C$. We obtain
a binary cubic homogeneous form in the two
remaining variables,  whose twisted discriminant $\Delta(U,V,W)$
is  the equation 
of the dual curve $\hat C$ (up to a square). This is because
this discriminant  cancels exactly 
when the intersection
$D.C$ has multiplicities.

\section{Pseudo-parameterizations}\label{sec:parame}

Let  $C$ be 
an  absolutely integral
plane projective curve over a field $k$. A {\it parameterization} of $C$
is a non-constant  map from $\PP^1$ onto $C$. In more concrete
terms  we have
a point $P_t=[X(t):Y(t):Z(t)]$ on $C$, depending on one 
formal parameter  $t$, the three projective coordinates beeing polynomials
in $k[t]$. It is well known \cite[theorem 4.11.]{SWP}
that a necessary condition
for such a parameterization  to  exist is that 
 $C$ has geometric genus zero. In particular this  never happens
for  an elliptic curve. One may relax the condition that the 
coordinates
$X(t)$, $Y(t)$ and $Z(t)$ should lye in $k[x]$ 
and allow for more general algebraic functions.
A typical restriction would be to ask that 
$X(t)$, $Y(t)$ and $Z(t)$ should  belong to 
a radicial extension of $k(t)$. In other words they should be
 rational fractions in $t$ and $\sqrt[e]{R(t)}$ for some positive
integer $e$ and some  $R(t)$ 
in $k(t)$. As explained in the introduction we will be interested in the
case when $C$ is a smooth cubic,   $k$ is a field  with characteristic 
different from $2$ and $3$, and $e=3$. We want to parametrize  plane
cubics by cubic radicals. Such a parameterization will be called
a {\it pseudo-parameterization} to avoid any confusion with
rational parameterizations that do not exist for genus one
curves.
We will assume that $C(k)$ is non-empty.
This is not a restriction if $k$ is a finite field.
We will even assume that $C$ has
a $k$-rational flex $O$.
This is not a restriction either, because every
cubic with a rational point is $k$-isomorphic to a plane cubic
with a rational flex.

We sketched  in the introduction how we claim to find pseudo-parameterizations.
We consider a line \[D_t : U(t)X+V(t)Y+W(t)Z=0\] 
in $\PP$, depending on one rational 
parameter $t$. 
Since every line in $\PP$ corresponds to a point in $\hPP$ we can 
associate to the family $D_t$ a rational curve $L\subset \hPP$
which is the image of the map
\begin{equation}\label{eq:paramL}
t\mapsto [U(t):V(t):W(t)].
\end{equation}

We saw in Section~\ref{section:dual} that the intersection
$D_t.C$ is described by a cubic form whose twisted discriminant
$\Delta(t)$ is, up to a square, equal to    $G(U(t),V(t),W(t))$ where
$G(U,V,W)=0$ is the projective equation of the dual $\hat C$.
So we look for polynomials $U(t)$, $V(t)$ and $W(t)$ such that 
$G(U(t),V(t),W(t))$ is a square in $k(t)$. A geometric
interpretation of the latter condition is that the rational curve
$L$ meets the dual $\hat C$ with all
even multiplicities. 
So we look for a rational curve $L\subset \hat\PP$
that intersects the dual curve $\hat C$
with  even multiplicities. Such a rational curve 
may be given by its
projective equation, or as the image of a parameterization as 
in~(\ref{eq:paramL}).

One may wonder if every pseudo-parameterization
occurs in that way. We briefly explain why this is
essentially the case. A pseudo-parameterization $t\mapsto P_t$
is a surjective map from a cyclic covering of $\PP^1$
onto $C$. So we have two conjugated points $P'_t$ and
$P''_t$. Since $C$ has a rational flex $O$, we have
 a chord and tangent group law, denoted $\oplus$, on it. We consider
the  sum $Q_t=P_t\oplus P'_t \oplus P''_t$. This is a point on $C$ defined over
$k(t)$, or equivalently  a map $t\mapsto Q_t$. We saw that such a 
map must be constant because $C$ has genus $1$.
So $P_t\oplus P'_t\oplus P''_t$ is a constant point $A\in C(k)$.
If $A$ is the origin $O$ then for every value of the parameter
$t$, the three points $P_t$, $P'_t$ and $P''_t$ are colinear.
They lye on a line  $D_t$ with equation $U(t)X+V(t)Y
+W(t)Z=0$ where $U(t)$, $V(t)$ and $W(t)$ are in
$k[t]$. So the pseudo-parameterization $t\mapsto P_t$
is of the type  studied above. 
If $A$ is not $O$,  we may look for  a point $B\in C(k)$ such that
$B\oplus B\oplus B=A$. Such a point always exists if $k$
is a finite field and $\# C(k)$ is not divisible by $3$.
Then   we  set
$R_t=P_t\ominus B$ and check that  $R_t\oplus R'_t\oplus R''_t=O$.
So the pseudo-parameterization $t\mapsto P_t$  is of the type
studied above, up to translation by a constant factor.
In general, we set  $R_t=P_t\oplus P_t\oplus P_t\ominus A$ and check
that $R_t+R'_t+R''_t=O$. So $t\mapsto P_t$ is of the type
studied above, up to a translation and a multiplication
by  $3$ isogeny.

We will say that two pseudo-parameterizations $t\mapsto P_t$
and $t\mapsto Q_t$ are {\it equivalent}
if there exists a birational fraction $\phi(t)$ 
such that $Q_{t}=P_{\phi(t)}$. 
We may  wonder if two different  families of projective lines
$t\mapsto D_t$ and $t\mapsto E_t$ can  give rise to equivalent
pseudo-parameterization  $t\mapsto P_t$ and $t\mapsto Q_t$.
In that case   $P_{\phi(t)}=Q_t$ 
lies in  the intersection of $D_{\phi(t)}$ and $E_t$.
If these two lines are distinct  then 
their intersection consists
of a single point  $P_{\phi(t)}=Q_t$  defined over $k(t)$.
Since every $k(t)$-rational point on $C$ is constant we deduce
that $P_t$ and $Q_t$ are constant. 
A contradiction. So $D_{\phi(t)}=E_t$ and the two families correspond
by a change of variable. In particular the two associated rational
curves in the dual plane are the same.

The conclusion is that finding pseudo-parameterizations boils down to
finding rational curves $L$ in the dual plane $\hPP$ having
even intersection with $\hat C$. It is natural to study first
rational curves going through several cusps of $\hat C$, because
the multiplicity intersection at  a singular point is 
greater than and generically equal to $2$. In the next section
we look for such rational curves with a low degree.

\section{The geometry of flexes}\label{sec:flexes}

Let $C\subset \PP$ be a smooth plane projective cubic. The nine flex points of
$C$ define a configuration in the plane $\PP$. More interestingly,
the nine flex tangents
correspond to nine points in the dual plane $\hat\PP$. We study the latter
configuration.
We are particularly 
interested in  low degree \emph{rational} curves going through
many of these nine cusps of $\hat C$. 
Remind  a \emph{rational} curve is  a  curve with  geometric genus
0 and a rational point. 
This is equivalent to the existence
of a rational parameterization, see \cite{SWP}, theorem 4.11. 
We will first assume that $C$ is the  Hessian plane cubic given
by  Equation~(\ref{eq:Hessian}).
Indeed, any smooth plane cubic can be mapped onto such an Hessian
cubic by a projective linear transform, possibly  after replacing $k$
by a finite extension of it. 
The modular invariant of $C$ is \[j(a)=\frac{27a^3(a+2)^3(a^2-2a+4)^3}{(a-1)^3(a^2+a+1)^3}.\]
%In particular $j\not =0$.
The nine flexes of $C$ are the three
points in the orbit of $O=(0:-1:1)$ under the action of $\cS_3$,
plus the six points in the orbit of $(-1:\zeta_3:0)$ 
under the action of $\cS_3$. 
Let \[\omega_{C} : (X:Y:Z)\mapsto (X^2-aYZ:Y^2-aXZ:Z^2-aXY)\]
be the Gauss map  associated with $C$.
The images by $\omega_{C}$ of the nine flexes 
are   the three
points in the orbit of  $(a:1:1)$ 
under the action of $\cS_3$ plus the 
six points in the orbit of $(\zeta_3^2:\zeta_3:a)$ 
under the action of $\cS_3$.  
Figure~\ref{fig:flexcusps}  lists these flexes
and their images by the Gauss map. We set $O=A_0=(0:-1:1)$ and $\hat O=B_0=
(a:1:1)$.

\begin{figure}[htb]
\centering
\begin{tabular}{|c|c|}
  \hline
  Flex of $C$ & Cusp on $\hat C$ \\
  \hline
  $A_0=(0:-1:1)$ & $B_0=(a:1:1)$  \\
  $A_1=(-1:1:0)$ & $B_1=(1:1:a)$  \\
  $A_2=(1:0:-1)$ & $B_2=(1:a:1)$  \\
  $A_3=(-1:\zeta_3:0)$ & $B_3=(\zeta_3^2:\zeta_3:a)$  \\
  $A_4=(\zeta_3:0:-1)$ & $B_4=(\zeta_3:a:\zeta_3^2)$  \\
  $A_5=(0:-1:\zeta_3)$ & $B_5=(a:\zeta_3^2:\zeta_3)$  \\
  $A_6=(\zeta_3:-1:0)$ & $B_6=(\zeta_3:\zeta_3^2:a)$  \\
  $A_7=(-1:0:\zeta_3)$ & $B_7=(\zeta_3^2:a:\zeta_3)$  \\
  $A_8=(0:\zeta_3:-1)$ & $B_8=(a:\zeta_3:\zeta_3^2)$  \\
  \hline
\end{tabular}
\caption{Flexes of $C$ and the corresponding cusps on its dual}
\label{fig:flexcusps}
\end{figure}

These nine points 
in the dual plane form an interesting configuration, depending
on the single parameter $a$. 
% One  can first  check,  e.g. using the Maple code given in
% Figure~\ref{fig:MapleI}, 

\paragraph{Position with respect to lines}
One  can first  check,  e.g. by exhaustive search, 
that no three among these  nine cusps  in the dual plane are 
colinear unless the modular invariant is zero. See the proof of 
Proposition~1 in Section~7.2 of \cite{BK}.
So the nine points in the dual plane corresponding to the nine flex lines are in general position
with respect to lines. 
We deduce the following lemma by duality.

\begin{lemma}
A smooth plane projective cubic over a field with prime to six
characteristic has no three concurrent
tangent flexes, unless its modular invariant
is zero.
\end{lemma}

\paragraph{Position with respect to conics} We now consider the configuration of the nine  flex tangents from the point of view of pencils
of conics. Remember that six points in general position do not lie on any conic.
Six pairwise distinct 
points lying on a conic are said to be {\it coconic}.
Six pairwise distinct lines  are said to be {\it coconic}
if they all are tangent to  a smooth conic.

\begin{lemma}\label{lem:conics}
Consider a smooth plane projective cubic over a field with prime to six
characteristic and assume that  its modular invariant
is not zero. Remove $3$ colinear flex points. The six  tangents  at
the six remaining flexes are coconic. There are twelve such 
configurations of six coconic flex tangents.
\end{lemma}

Note that we
claim that the six flex tangents are coconic. Not the six flex points.
Equivalently we claim that the six points in the dual plane corresponding to the six
flex tangents  are coconic.

We  first note that  the conic with equation $UW-aV^2=0$ meets
$\hat C$ at $(a:1:1)$, $(1:1:a)$, $(\zeta_3^2:\zeta_3:a)$,  
$(a:\zeta_3^2:\zeta_3)$, $(\zeta_3:\zeta_3^2:a)$,  and
$(a:\zeta_3:\zeta_3^2)$. The three remaining flexes in
$\PP$ are $(1:0:-1)$, $(\zeta_3:0:-1)$ and
$(\-1:0:\zeta_3)$ and they lie on the line with equation
$Y=0$.
The  action of $\cS_3$ produces two more similar conics.

The conic with equation $U^2+V^2+W^2+(a+1)(UV+UW+VW)=0$
meets $\hat C$ at the six points in the
orbit of $(\zeta_3^2:\zeta_3:a)$  under the action of $\cS_3$.
The three remaining flexes in
$\PP$ are $(0:-1:1)$, $(-1:1:0)$, and $(1:0:-1)$.
They lie on the line with equation $X+Y+Z=0$.

The conic with equation $U^2+\zeta_3V^2+\zeta_3^2W^2
+(a+1)(\zeta_3^2UV+\zeta_3UW+VW)=0$
meets $\hat C$ at the three  points in the
orbit of $(a:1:1)$  under the action of $\cS_3$. And also
at the three  points in the
orbit of $(\zeta_3^2:\zeta_3:a)$  under the action of $\cS_3$. 
The three remaining flexes in
$\PP$ are $(0:\zeta_3:-1)$, $(\zeta_3:-1:0)$, and $(-1:0:\zeta_3)$.
They lie on the line with equation $X+\zeta_3Y+\zeta_3^2Z=0$.
The  action of $\cS_3$ produces one more such conic.

The conic with equation $\zeta_3U^2+V^2+\zeta_3W^2
+(a+\zeta_3^2)(UV+\zeta_3^2UW+VW)=0$
meets $\hat C$ at 
$(a:1:1)$, $(1:1:a)$, $(\zeta_3:a:\zeta_3^2)$,
$(a:\zeta_3^2:\zeta_3)$, $(\zeta_3:\zeta_3^2,a)$, 
$(\zeta_3^2:a:\zeta_3)$.
The three remaining flexes in
$\PP$ are $(1:0:-1)$, $(-1:\zeta_3:0)$, and $(0:\zeta_3:-1)$.
They lie on the line with equation $\zeta_3X+Y+\zeta_3Z=0$.
The  action of $\cS_3$ produces five more conics.

We thus obtain twelve smooth 
conics that cross the dual curve $\hat C$ at
six out of its nine cusps. Each of these conics is associated
with one of the twelve triples of colinear flexes. \hfill $\Box$

Four among  these twelve 
conics are especially
interesting because their equations do not involve $\zeta_3$.
We note that  three  among these four conics are clearly
rational over $k(a)$ because they have an evident $k(a)$ rational point. 
The last one is rational also because its quotient by the evident
automorphism of order $3$ is $\PP^1$ over $k(a)$.

\paragraph{Position with respect to cubics}
Next we study the pencil of cubics going through the nine  points 
in the dual plane associated with
the nine flex tangents. It has projective dimension zero in general.
The cubic with equation
\[a(U^3+V^3+W^3)=(a^3+2)UVW\]
goes  through all these nine points in the dual
plane. 
This cubic is  in  general non-singular. So it is not
particularly interesting for our purpose. 
%The corresponding   Maple code is given in  Figure~\ref{fig:MapleIII}.

\paragraph{Position with respect to quartics} We now consider curves of degree $4$ in the dual plane. The projective
dimension of the space of plane quartics is $14$. So we can force
a quartic to meet the  $9$ points we are interested in and there
remains $5$ degrees of freedom. Since we are particularly interested
in rational  curves we use these remaining degrees of freedom
to impose a big singularity at $\hat O = B_0=(a:1:1)$. Indeed, two
degrees of freedom suffice  to  cancel the degree
$1$ part in the Taylor expansion at $\hat O$. And three more
degrees of freedom suffice  to  cancel the degree $2$
part also. We find a rational quartic $Q$ in $\hat \PP$
passing through the nine cusps of $\hat C$
and having intersection multiplicity at least 
two at each  of them (because they
are cusps) and  at least six at the cusp $\hat O$. 
%This latter multiplicity is actually
%equal to  four.  Indeed the intersection $Q.C$ consists of
%the nine flexes  plus a degree $3$ divisor. Since the
%flexes lie on a cubic distinct from $C$ (e.g. any  cubic in the
%Hessian family), the remaining $3$ points must form the
%intersection of $C$ with  a line. Since two among  these points are
%already equal to $(a:1:1)$, the line in question must be the (flex) tangent
%at $(a:1:1)$.
The equation
of this rational quartic $Q$ is

\begin{equation*}
\begin{split}
 U^4+a(V^4+W^4)-2a(U^3V+U^3W+V^3W+VW^3)-(a^3+1)U(V^3+W^3)\\
+3a^2U^2(V^2+W^2)+(a^4+2a)V^2W^2+(1-a^3)UVW(V+W)=0.
\end{split}
\end{equation*}

This quartic
is irreducible  as soon as the modular invariant of
$C$ is non-zero, which we assume from now on. 
Computing the intersection with all lines
through $\hat O$ we find the following 
 parameterization of this  quartic 

\begin{eqnarray*}
U(t)&=& a^2t^4-2at^3+(a^3+2)t^2-2a^2t+a,\\
V(t)&=& a^4t^4+(1-3a^3)t^3+3a^2t^2-2at+1,\\
W(t)&=&  at^4-(a^3+1)t^3+3a^2t^2-2at+1.
\end{eqnarray*}

\begin{figure}[!h]
\centering
{\includegraphics[width=4cm,height=4cm]{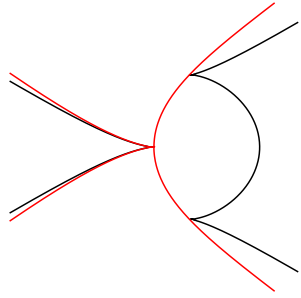}}
\caption{The real part of the intersection of $\hat C$ and $Q$.}
\label{fig:qc}
\end{figure}

Substituting $U$, $V$, and $W$  by $U(t)$, $V(t)$, and 
$W(t)$ in the equation of $\hat C$ we find the  degree 
$24$ polynomial 
\[t^6(t+1)^2(t^2-t+1)^2(at-2)^2((a+1)t-1)^2((a^2-a+1)t^2+(1-2a)t+1)^2
(a^2t^2+1-at)^2.\]
We check that
 $Q$ has two branches at $\hat O$.
One branch corresponds to $t=0$, and it has intersection
multiplicity  $6$ with $\hat C$.
The other branch corresponds to  $t=2/a$, and it has  
intersection multiplicity
$2$ with $\hat C$. This is illustrated by Figure~\ref{fig:qc} where 
the real locus of $\hat C$ is in 
black and the real locus of $Q$ is in red. So the total multiplicity of $Q.\hat C$ at
$\hat O$ is
$8$. And the intersection $Q.\hat C$ only consists of cusps 
of $\hat C$; one with multiplicity $8$ and the eight 
others with multiplicity $2$.
The real part of this intersection
locus is visible on  Figure~\ref{fig:qc}.
%The corresponding   Maple code is given in  Figure~\ref{fig:MapleIV}.

\begin{lemma}\label{lem:quartic}
Consider a smooth plane projective cubic $C$ over a field with prime to six
characteristic and assume  that its modular invariant
is not zero. 
Let $\hat C$ be the dual of $C$.
Let $\hat O$ be one of the nine cusps of $\hat C$.
There exists a rational quartic $Q$ in the dual plane, 
such that the intersection
$Q.\hat C$ has multiplicity $8$ at $\hat O$ and 
$2$ at each of the eight remaining cusps.
In particular $Q.\hat C$ is an even combination of
cusps of $\hat C$.
\end{lemma}

We stress that the definition of the quartic $Q$
involves one flex
on the one hand, and the eight remaining flexes on the other hand.
So we can define this quartic
 for any cubic having a rational  flex,
that is for any elliptic curve (and this makes a difference with
the four conics constructed earlier, that distinguish
a triple of  colinear
flexes, and therefore cannot always 
be defined over the base field.)

So we  can take  for $C$ an elliptic curve with Weierstrass equation
\begin{equation}\label{eq:weier}
F(X,Y,Z)=Y^2Z-X^3-aXZ^2-bZ^3.
\end{equation}
We assume  $a\not =0$, so the modular invariant
is non-zero either. The image of the origin $O=(0:1:0)$
by the Gauss map is $\hat O = (0:0:1)$,
and the  quartic $\hat Q$ given by Lemma~\ref{lem:quartic} has equation
\[U^4-3V^4+6UV^2W=0,\]
and parameterization

\begin{eqnarray}\label{eq:parat}
U(t)&=& 6t^2,\\\nonumber
V(t)&=& 6t^3,\\\nonumber
W(t)&=& 3at^4-1.
\end{eqnarray}

\section{Intersecting a cubic with lines}\label{sec:inter}

In this section we assume that the map $a\mapsto a^3$ from
$k$ to $k$ is surjective. 
This is the case if $k$ is the field
of real numbers for example.  This is also the case if $k$
is a finite field with $q$ elements when $q$ is congruent to
$2$ modulo $3$.
For every element $a$ in $k$ we choose once and for all
a cubic root $\sqrt[3]{a}$ of $a$. This way we define a map
$\sqrt[3]{} : k\rightarrow k$. 
% Let $\zeta_3\subset \bar k$ be  a primitive third root of unity. %We set $l=k(\zeta_3)$.
We will use the general recipe in Section~\ref{sec:parame} and the rational
curves exhibited in Section~\ref{sec:flexes} to produce several pseudo-parameterizations
of a plane cubic $C$.

\subsection{Intersecting the dual curve with a conic}

We may first take $L$ to be one of the twelve conics in
Lemma~\ref{lem:conics}.
So we  assume that $C$ is the Hessian cubic given by
Equation~(\ref{eq:Hessian}) for some $a$ such that  $a^3\not =1$.
Four conics, among the twelve conics
given in
Lemma~\ref{lem:conics},
 are  rational over $k(a)$.
The intersection $L.\hat C$ has degree $12$ and contains
six among the nine cusps of $\hat C$, each with
multiplicity  $2$.
So this intersection is exactly twice the sum of these
six cusps.
If we take for $L$ the conic with equation
$UW-aV^2=0$ then a convenient parameterization is given
by $U(t)=1$, $V(t)=-t$ and $W(t)=at^2$. The corresponding
line $D_t$ has equation \[X-tY+at^2Z=0.\]
We substitute $X$ by $tY-at^2Z$ in the Hessian
Equation~(\ref{eq:Hessian}) and find the degree 
$3$ form in $Y$ and $Z$
\[(t^3+1)Y^3-3at(t^3+1)Y^2Z+3a^2t^2(t^3+1)YZ^2+(1-a^3t^6)Z^3\]
describing the intersection $C.D_t$. We  
divide by $(t^3+1)Z^3$ and we obtain a cubic polynomial in $y=Y/Z$ 
whose
twisted discriminant is \[\Delta(t)=\left(\frac{9(1+a^3t^3)}{1+t^3}  \right)^2.\]

We use the formulae and notation
in Section~\ref{sec:cubiceq}. We have
\begin{eqnarray*}
s_1&=&{3at},\\
s_2&=&3a^2t^2,\\
s_3&=&\frac{a^3t^6-1}{t^3+1},\\
\delta&=&\frac{9(1+a^3t^3)}{1+t^3} ,\\
R&=&-27\frac{a^3t^3+1}{t^3+1},\\
R'&=&0.
\end{eqnarray*}

So we find the solution
\[y=at-\sqrt[3]{\frac{a^3t^3+1}{t^3+1}},\]
and we deduce
\[x=X/Z=ty-at^2=-t\sqrt[3]{\frac{a^3t^3+1}{t^3+1}},\]

This is the  
pseudo-parameterization found by Farashahi~\cite{fara}.

\subsection{Intersecting the dual curve with a quartic}

Assume now that we  take $L$ to be the rational 
quartic   $Q$ in Lemma~\ref{lem:quartic}.
All the multiplicities
in the intersection $Q.\hat C$ are
even. So we expect the twisted discriminant  to be  a 
square. 
This time we may as well take for $C$ the Weierstrass cubic in
Equation~(\ref{eq:weier}). 
The 
parameterization of $Q$ given in Equation~(\ref{eq:parat}) provides
a  one parameter family  of 
lines $(D_t)_t$ with equation
\[6t^2X+6t^3Y+(3at^4-1)Z=0.\] We divide by $Z$, we set 
$x=X/Z$, $y=Y/Z$ and we substitute  $y$ by $1/(6t^{3})-at/2 -x/t$ in 
the Weierstrass
Equation~(\ref{eq:weier}). We
 find a cubic equation $x^3-s_1x^2+s_2x-s_3$ in $x=X/Z$, where
\begin{eqnarray*}
s_1&=&{1/t^2},\\
s_2&=&1/(3t^4),\\
s_3&=&(1/t^{6}-6a/t^2-36b+9a^2t^2)/36.\\
\end{eqnarray*}
Using the formulae and notation
in Section~\ref{sec:cubiceq} we find
\begin{eqnarray*}
\delta&=& (-1/t^6-108b-18a/t^2+27a^2t^2)/12,\\
R&=&0,\\
R'&=&(-1/t^6-108b-18a/t^2+27a^2t^2)/4.
\end{eqnarray*}

So we find the solution
\[x=X/Z=\frac{1}{3t^2}+\sqrt[3]{\frac{a^2t^2}{4} -\frac{1}{108t^6} 
-b-\frac{a}{6t^2} }\]
and 
\[y=Y/Z= \frac{1}{6t^{3}}-at/2 -x/t.\]

This is the  
pseudo-parameterization found by Icart~\cite{icart}, up
to the change of variable $t \leftarrow -1/t$. 

\subsection{Intersecting the dual curve with a line}

\begin{figure}[htb]
\centering
{\includegraphics[width=4cm,height=4cm]{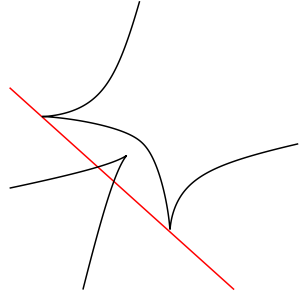}}
\caption{The intersection of $\hat C$ and $L$}
\label{fig:3tor}
\end{figure}

Assume finally that we take for $L$ a 
line passing through  two rational cusps
of $\hat C$. 
So  we  assume that $C$ is the Hessian
cubic given by Equation~(\ref{eq:Hessian}) for some
 $a^3\not =1$.
Assume  $L$ is the unique line passing through
the two cusps $B_0=(a:1:1)$ and $B_2=(1:a:1)$ of $\hat C$.
The intersection $L.\hat C$ has degree $6$.
Since  $(a:1:1)$ and $(1:a:1)$   each have intersection multiplicity $\ge 2$,
there  remains at most two intersection points. 
An illustration of this situation in the real projective
plane is given on
Figure~\ref{fig:3tor}.

Not all the multiplicities in the intersection $L.\hat C$
are even, but only two multiplicities are odd. So 
we expect $\Delta(t)$ to be a  square times a degree $2$ polynomial in
$t$. 
Points on $L\subset \hat \PP$ represent a linear  pencil of
lines in $\PP$ generated by the tangents to $C$ at
$(0:-1:1)$ and $(1:0:-1)$. 
The first tangent has equation $aX+Y+Z=0$. 
The second tangent has equation $X+aY+Z=0$. 
So let $t$ be a formal parameter and consider the line
$D_t$ with equation $(at+1)X+(t+a)Y+(t+1)Z=0$.
The tangent at $(0:-1:1)$ corresponds to the value $t=\infty$.
The tangent at $(1:0:-1)$ corresponds to the value $t=0$.
The line $D_t$ meets the fixed  point $(1:1:-a-1)$ and the moving
point 
$(1,-t,t-1)$. So a  parametric description  of $D_t$ is given
by \[i\mapsto (i+1:i-t:t-1-(a+1)i).\]
We substitute $X$ by $i+1$, $Y$ by $i-t$ and $Z$
by $t-1-(a+1)i$ in Equation~(\ref{eq:Hessian}) and divide by
the leading coefficient. We find the degree three polynomial
\begin{equation}\label{eq:hi}
h(i)=i^3+\frac{3t(a+2)i}{a^2+a+1}+\frac{3t(1-t)}{a^2+a+1}
\end{equation}
defining the intersection $D_t.C$. The twisted discriminant of 
$h$ is 
\begin{equation}\label{eq:delt}
\Delta(t)=81t^2\frac{9(a^2+a+1)t^2+2(2a+1)(a^2+a+7)t+9(a^2+a+1)}{(a^2+a+1)^3}.
\end{equation}

This is not quite a square in $k(a)(t)$.
However, it only has   two roots
with  odd multiplicity. So if we
  substitute $t$ by a well chosen rational fraction, 
we can
 turn $\Delta$ into a
square. 
So we look for a parameterization of the plane projective
conic  with equation
\begin{equation}\label{eq:conic}
(a^2+a+1)S^2=9(a^2+a+1)T^2+2(2a+1)(a^2+a+7)TK+9(a^2+a+1)K^2.
\end{equation}

This conic has two evident $k$-rational points, namely $(3:1:0)$ and
$(3:0:1)$. The line through these two points has equation
\[-S+3T+3K=0.\]
The tangent at $(3:0:1)$ has equation
\[3(a^2+a+1)S-(2a+1)(a^2+a+7)T-9(a^2+a+1)K=0.\]
The generic line in the linear pencil generated by these two lines has
equation
\begin{equation}\label{eq:genlin}
(3(a^2+a+1)-j)S+(3j-(2a+1)(a^2+a+7)j)T+(3j-9(a^2+a+1)j)K=0
\end{equation}
where $j$ is a formal parameter.

Intersecting the conic in Equation~(\ref{eq:conic}) with the line
in Equation~(\ref{eq:genlin}) we find the parameterization

\[\left\{
\begin{array}{ccc}
S(j) &=&3j^2-2(a+2)^3j+3(a+2)^3(a^2+a+1),\\
T(j) &=& j(j-3(a^2+a+1)), \\
K(j) &=&(a^2+a+1)((a+2)^3-3j).
\end{array}
\right.\]

We now substitute $t$ by $T(j)/K(j)$ in Equation~(\ref{eq:hi})
and find a cubic polynomial with coefficients in the field $k(a)(j)$.
If we  substitute $t$ by $T(j)/K(j)$ in Equation~(\ref{eq:delt})
we  find that $\Delta =\delta ^2(j)$
where
\[\delta(j)=
\frac{9j(3j^2-2(a+2)^3j+3(a^2+a+1)(a+2)^3)(3(a^2+a+1)-j)}{((a+2)^3-3j)^2(a^2+a+1)^3}.\]

We use the formulae and notation
in Section~\ref{sec:cubiceq}. 
The polynomial $h$ in Equation~(\ref{eq:hi}) has
coefficients $1$, $-s_1$, $s_2$ and $-s_3$
with \begin{eqnarray*}
s_1&=&0\\
s_2&=&-\frac{3j(a+2)(3(a^2+a+1)-j)}{(a^2+a+1)^2((a+2)^3-3j)}\\
s_3&=&\frac{3j(3(a^2+a+1)-j)((a^2+a+1)(a+2)^3-j^2)}{(a^2+a+1)^3((a+2)^3-3j)^2}.
\end{eqnarray*}

We deduce the following
pseudo-parameterization   of the cubic $C$

\begin{eqnarray*}
R(j)&=&\frac{27j^2(3(a^2+a+1)-j)}{((a+2)^3-3j)(a^2+a+1)^3}\\
\rho(j)&=&\sqrt[3]{R(j)}\\
\rho'(j)&=&\frac{9j(a+2)(3(a^2+a+1)-j)}{(a^2+a+1)^2((a+2)^3-3j)\rho(j)}\\
i(j)&=&\frac{\rho(j)+\rho'(j)}{3}\\
t(j)&=&\frac{j(3(a^2+a+1)-j)}{(a^2+a+1)((a+2)^3-3j)}\\
P(j)&=&(i(j)+1:i(j)-t(j):t(j)-1-(a+1)i(j)).
\end{eqnarray*} 
where $P(j)$ 
is the point
on $C$ associated with the parameter $j$.

\begin{figure}[hh!]
\centering
{\includegraphics[width=4cm,height=4cm]{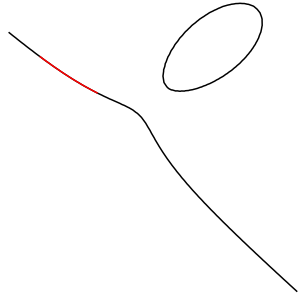}}
\caption{A pseudo-parameterization}
\label{fig:paramcub}
\end{figure}

We illustrate this situation on Figure~\ref{fig:paramcub}
in the case $a=2$.
The red  segment corresponds to the parameter $j$
taking  values in the interval $[-4,-0.3]$. 
We also note that  the computation in Section~3.1 of
\cite{KLR} hides a similar geometric situation.

\section{Classifying pseudo-parameterization}\label{sec:K3}

We have seen many different pseudo-parameterizations of
a plane cubic, each associated
with a rational curve in $\hPP$  having even intersection with
the dual curve $\hat C$  in Equation~(\ref{eq:G}).
We may wonder if there exist more such rational curves,
 leading to more pseudo-parameterizations. We may also try to
put some structure on the set of such curves. This is
our  purpose of this section. We assume that the reader has some familiarity
with algebraic surfaces as presented in \cite{Z,bad}, and particularly
with elliptic  and K3 surfaces \cite{ST,elk,bertin}. We shall not enter
into the details.
Any  rational curve $L$ having even intersection
with $\hat C$ lifts to a rational curve
on the degree two covering $\Sigma$ of $\hPP$
branched along  $\hat C$. To define $\Sigma$ we consider
the function field $k(a)(U/W,V/W)$ of $\hPP$ over $k(a)$. We
define a quadratic extension of this field by adding a 
square root $\gamma$ of $G(U,V,W)/W^6$ where $G(U,V,W)$
is the equation of $\hat C$. The normal
closure of $\hPP$ inside $k(a)(U/W,V/W,\gamma)$  is $\Sigma$.
It has nine singularities. One
above each of the nine cusps of $\hat C$. 
In order to obtain a smooth model for $\Sigma$, we first
blow up $\hPP$ at each of the cusps of $\hat C$. We 
call $\Pi$ the resulting surface. The inverse image of
$\hat C$ by  $\Pi\rightarrow \hPP$ consists of one  smooth genus one curve 
and $9$ rational curves tangent to it.
We call $S$ the normal closure of $\Pi$ in  $k(a)(U/W,V/W,\gamma)$.
This is a smooth surface, the minimal model of $\Sigma$.

We call $\sigma_1$ the automorphism of 
$\hPP$ that maps $[U:V:W]$ onto $[V,W,U]$.
We call $\sigma_2$ the automorphism of 
$\hPP$ that maps $[U:V:W]$ onto $[U,\zeta_3 V,\zeta_3^2W]$.
We call $\sigma_3$ the automorphism of 
$\hPP$ that maps $[U:V:W]$ onto $[V,U,W]$.
We extend these three automorphisms
to $k(a)(U/W,V/W,\gamma)$ by sending $\gamma$
to itself. The resulting automorphisms are  called $\sigma_1$,
$\sigma_2$ and $\sigma_3$ also.
They  induce automorphisms of $\Pi$, $\Sigma$
and $S$ denoted $\sigma_1$,
$\sigma_2$ and $\sigma_3$  again.
We call $\sigma_4$ the unique non-trivial automorphism 
of $k(a)(U/W,V/W,\gamma)$ over $k(a)(U/W,V/W)$. It induces
automorphisms of $\Sigma$ and $S$ denoted $\sigma_4$.
The action of $\sigma_1$,
$\sigma_2$, $\sigma_3$ on the $B_i$  is given
by the following three permutations  of the indices
\begin{eqnarray*}
\sigma_1&=&(0,1,2)(3,4,5)(6,7,8),\\
\sigma_2&=&(0,5,8)(1,3,6)(2,4,7),\\
\sigma_3&=&(0,2)(1)(3,6)(4,8)(5,7).
\end{eqnarray*}

The group generated by $\sigma_1$ and $\sigma_2$ has
order nine. It acts simply transitively on the nine cusps, and
also on the nine corresponding rational curves  on the blow up 
$\Pi$. We choose one of the two rational  curves on $S$
above $B_0$ and call  it $E_0$. For $1\le i\le 9$
we call $E_i$ the image of $E_0$ by the unique automorphism
in $<\sigma_1,\sigma_2>$ that maps $B_0$ onto $B_i$.
We call $F_i$ the image of $E_i$ by $\sigma_4$.
We thus obtain eighteen rational curves on $S$. Let $H$ be the inverse
image by $S\rightarrow \hPP$ of any line in $\hPP$. The lattice
generated by the $E_i$, $F_i$ and $H$ in the Néron-Severi group
has  rank $19$, and discriminant $2.3^9$. The intersection
indices are 
\begin{eqnarray*}
E_i.F_i&=&1,\\
E_i^2&=&-2,\\
F_i^2&=&-2,\\
E_i.E_j&=&0 \text{ for } i\not =j,\\
E_i.F_j&=&0 \text{ for } i\not =j,\\
E_i.H&=&0,\\
F_i.H&=&0,\\
H^2&=&2.
\end{eqnarray*}

Let $D$ be a generic
line in $\hPP$ through $B_0$.
The intersection of $D.\hat C$ is $2B_0$ plus
an effective degree four divisor.
So the inverse image of $D$ in $S$ is the union of $E_0$, $F_0$
and a genus one curve
with at least two rational points : the  intersection points
with $E_0$ and $F_0$.
Thus the inverse image by $S\rightarrow \hPP$ 
of the pencil of  lines through $B_0$ defines  an elliptic
fibration $f : S\rightarrow \PP^1$ of $S$, with  two  sections $E_0$ and $F_0$,
so  $S$ is an elliptic K3 surface. 
The following lemma \cite[2.3]{these} is usefull when
looking for rational curves on a K3 surface.

\begin{lemma}
Let $D$ be a class with self-intersection $-2$ in the Néron-Severi group of a K3 surface.
Then either $D$ or $-D$ contains an  effective divisor.
If this divisor is irreducible then it is a smooth rational
curve.
\end{lemma}

We may also look for singular rational curves in classes with positive
self-intersection. One can even count rational
curves in such classes \cite{beau,Lee,Wu}. Since there are many of them,
they are unlikely  to be defined over the base field.
Indeed, all the rational curves in Section~\ref{sec:flexes} lift
to smooth rational curves on $S$ having self-intersection $-2$.
For example the conic in $\hPP$ passing through $B_0$, $B_1$, $B_2$, $B_3$, $B_4$, $B_5$ lifts to a rational
curve $I_{012345}$ on $S$. We have 
$H.I_{012345}=2$, $E_0.I_{012345}=E_1.I_{012345}=E_2.I_{012345}=1$
and  $F_3.I_{012345}=F_4.I_{012345}=F_5.I_{012345}=1$ and
$I_{012345}$ has zero  intersection with the remaining
$E_i$ and $F_i$.
We deduce  the following identity 
in the Néron-Severi group
\[3I_{0,1,2,3,4,5}=3H-2(E_0+E_1+E_2)-(F_1+F_2+F_3)
-(E_3+E_4+E_5)-2(F_3+F_4+F_5),\]
and $I_{0,1,2,3,4,5}$ has  self-intersection $-2$.
We find similarly, and with evident notation,
\[3I_{0,1,3,4,7,8}=3H-2(E_0+E_3+E_7)-(F_0+F_3+F_7)
-(E_1+E_4+E_8)-2(F_1+F_4+F_8),\] and
\[3I_{0,1,3,5,6,8}=3H-2(E_0+E_5+E_8)-(F_0+F_5+F_8)
-(E_1+E_3+E_6)-2(F_1+F_3+F_6).\]
The action
of  $<\sigma_1,\sigma_2,\sigma_3, \sigma_4>$ produces
$24$ similar smooth rational curves on $S$
with self intersection
$-2$. This is the contribution  of conics in Lemma~\ref{lem:conics}.

Now consider the quartic given by Lemma~\ref{lem:quartic}. It lifts
to a rational curve $J_0$ on $S$, such  that $J_0.H=4$,
$J_0.E_0=2$, $J_0.F_0=1$,
$J_0.E_i=1$, $J_0.F_i=0$  for $1\le i\le 8$. We have
the following identity in the  Néron-Severi group
\[3J_0=6H-5E_0-4F_0-\sum_{1\le i\le 8}(2E_i+F_i).\]

\begin{figure}[!h]
\centering
{\includegraphics[width=6cm,height=6cm]{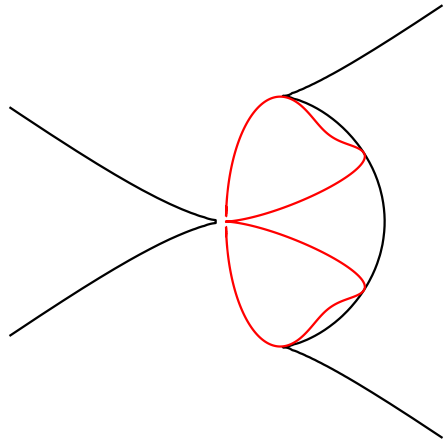}}
\caption{One more rational curve having  even intersection with $\hat C$.}
\label{fig:nouv}
\end{figure}

The action of  $<\sigma_1,\sigma_2,\sigma_4>$ produces
$18$ such rational curves with self intersection $-2$. The lattice generated
by $H$, the nine   $E_i$, the nine  $F_i$, and 
the $24+18$ classes coming from conics and quartics, has dimension
$19$ and discriminant $54$. This is the full
Néron-Severi group  of $S$ when $k$ has characteristic zero and
$a$ is a transcendental.
Using the knowledge of this  Néron-Severi group we can prove  that
there are infinitely many rational curves on $S$, leading to infinitely
many pseudo-parameterizations of the cubic $C$.
We consider an elliptic-fibration of $S$, for example the fibration
$f : S \rightarrow \PP^1$ introduced above. We choose the section
$E_0$ as origin.
The generic fiber of $f$
is an elliptic curve over the function field $k(t)$ of $\PP^1$.
Fibers of $f$ map onto lines through $B_0$ in $\hPP$. The
height singular fibers of $f$ map onto the lines $B_0B_i$
for $1\le i\le 8$. Each of them has Kodaira type $I_3$, the three irreducible 
components being $E_i$, $F_i$, and a third rational curve $G_i$ crossing
$E_0$ and $F_0$.
Let $T\subset \NS(S)$ 
be the group generated by the zero section $E_0$ and the fiber
components $E_i$, $F_i$, $G_i$ for $1\le i\le 8$.
The Mordell-Weil group of the generic fiber is isomorphic \cite[Theorem 6.3]{ST}
 to the quotient
$\NS(S)/T$. Since $E_i+F_i+G_i=H-E_0-F_0$ does not depend on $i$ for $1\le i\le 8$,
the rank of $T$ is $18$ and the rank of  $\NS(S)$ is one.
So we have infinitely many sections of $f$. The images of these
sections all are  rational curves with self  intersection $-2$.
We draw one of these rational curves (rather its image in $\hPP$)
on Figure~\ref{fig:nouv}. In case $C$ is the Weierstrass
cubic in Equation~(\ref{eq:weier}), a 
parameterization of this
rational curve is

\begin{eqnarray}%\label{eq:parat}
U(t)&=&4at^6 + 4t^2/27,\\\nonumber
V(t)&=& t(4at^6 + 4t^2/27),\\\nonumber
W(t)&=& a^2t^8 + 2at^4/27 + 4bt^6 + 1/81.
\end{eqnarray}

%Since we already found $19$ linearly independant classes in the Néron-Severi group, we know that the rank of $S$ is at least $19$.  The rank of an algebraic  K3 surface  is at most $20$ and this maximal is reached only for the Kummer surface of the product of two isogenous CM elliptic curves. This is not the case for $S$ because it is defined over $k(a)$ where $a$ is a transcendental. So the rank of 

\end{document}